\documentclass{amsart}

\usepackage{latexsym}
\usepackage{amssymb,amsmath,amsfonts}

\newtheorem{theorem}{Theorem}[section]

\theoremstyle{definition}

\newtheorem{example}[theorem]{Example}

\theoremstyle{remark}
\newtheorem{rmk}[theorem]{Remark}

\numberwithin{equation}{section}

\numberwithin{equation}{section}

\begin{document}

\title[Erratum]
 {A note on the paper: On iterations for families of asymptotically pseudocontractive mappings.}
 \address{$^1$ Department of Mathematics, Nnamdi Azikiwe University, P. M. B. 5025, Awka, Anambra State, Nigeria.}

 \email{ euofoedu@yahoo.com}
\author[E. U. Ofoedu]{ Eric U. Ofoedu$^{*\;1}$}

\keywords{Asymptotically pseudocontractive mappings, uniformly $L$-Lipschitzian mapping, real Banach space.\\{\indent 2010 {\it Mathematics Subject Classification}. 47H06,
47H09, 47J05, 47J25.\\$^*$The author undertook this work when he visited
 the Abdus Salam International Centre
for Theoretical Physics (ICTP), Trieste, Italy as a visiting fellow.}}

\begin{abstract}
It is our aim in this note to give a counter example to an argument used in the proof of the main theorem of the paper: On iterations for families of asymptotically pseudocontractive mappings, \emph{Applied Mathematics Letters,} {\bf 24} (2011), 33-38 by A. Rafiq \cite{Rafiq}; and give an alternative condition to correct the anomaly.
\end{abstract}

\maketitle
\section{Introduction.} 

\noindent This work is motivated by the recent paper of A. Rafiq \cite{Rafiq}. Careful reading of Rafiq's work shows that there is a serious gap in the proof of Theorem 5 of \cite{Rafiq}, which happens to be main theorem of the paper.\\

It is our aim to give a counter example to the argument used in the proof of Theorem 5 of \cite{Rafiq} and suggest an alternative condition in order to close the observed gap.

\section{Preliminary.}

\noindent Let $E$ be a real Banach space with dual $E^*$ and let $\big< ., \;.\big>$ be the duality pairing between members of $E$ and $E^*$. The mapping $J: E\rightarrow 2^{E^*}$ defined by
\begin{equation*}
J(x) = \{f^*\in E^* : \big<x,f^*\big> =||x||^2; ||f^*||=||x||\}, x\in E,
\end{equation*}
is called the normalized duality mapping. We note that in a Hilbert space $H$, $J$ is the identity operator. The single valued normalized duality mapping is denoted by $j$.\\

\noindent A mapping $T:D(T)\subset E\to E$ is said to be \emph{$L$-Lipschitzian} if there exists $L>0$ such that $$\|Tx-Ty\|\le L\|x-y\|\;\forall\;x,y\in D(T);$$ and $T$ is said to be \emph{uniformly $L$-Lipschitzian} if there exists $L>0$ such that 
$$\|T^nx-T^ny\|\le L\|x-y\|\;\forall\;x,y\in D(T),\forall\;\;n\ge 1,$$ where $D(T)$ denotes the domain of $T$. It is well known that the class of uniformly $L$-Lipschitzian mappings is \emph{a proper subclass} of the class of $L$-Lipschitzian mappings.\\

\noindent The mapping $T$ is said to be \emph{asymptotically pseudocontractive} if there exists a sequence $\{k_n\}_{n\ge 1}\subset [1,+\infty)$ with $\displaystyle \lim_{n\to\infty}k_n=1$ and for all $x,y\in D(T),$ there exists $j(x-y)\in J(x-y)$ such that 
$$\Big<T^nx-T^ny,j(x-y)\Big>\le k_n\|x-y\|^2\;\forall\;x,y\in D(T),\;\forall\;n\ge 1.$$

\noindent In \cite{Rafiq}, A. Rafiq studied the strong convergence of the sequence $\{x_n\}_{n\ge 1}$ defined by 

\begin{eqnarray*}
&&x_1\in K,\\
x_{n+1}&=&(1-\alpha_n)x_n+\alpha_nT_1^ny_n^1\\
y_n^i&=&(1-\beta_n^i)x_n+\beta_n^iT_{i+1}^ny_n^{i+1}\\
&\vdots&\\
y_n^{p-1}&=&(1-\beta_n^{p-1})x_n+\beta_n^{p-1}T_p^nx_n, \;n\ge 1,\;\;\;\;\;\;\;\;\;\;\;\;\;(1.4)
\end{eqnarray*} for approximation of common fixed point of finite family of asymptotically pseudocontractive mappings in real Banach space. He proved the following theorem.

\begin{theorem}(See Theorem 5 of \cite{Rafiq}) Let $K$ be a nonempty closed convex subset of a real Banach space $E$ and $T_l:K\to K,\;l=1,2,...,p;\;p\ge 2$ be $p$ asymptotically pseudocontractive mappings with $T_1$ and $T_2$ having bounded ranges and a sequence $\{k_n\}_{n\ge 1}\subset [1,+\infty),$ $\displaystyle\lim_{n\to\infty}k_n=1$ such that $x^*\in \displaystyle \bigcap\limits_{l=1}^{p}F(T_l)=\{x\in K:T_1x=x=T_2x=...=T_px\}.$ Further, let $T_1$ be uniformly continuous and $\{\alpha_n\}_{n\ge 1}$, $\{\beta_n^i\}_{n\ge 1},$ $\{\beta_n^{p-1}\}_{n\ge 1}$ be sequences in $[0,1], \;i=1,2,...,p;\;p\ge 2$ such that\\
 $(i)\;\displaystyle\lim_{n\to\infty}\alpha_n=0=\lim_{n\to\infty}\beta_n^1;$\\
 $(ii)\;\displaystyle \sum_{n\ge 1}\alpha_n=\infty.$ \\
 For arbitrary $x_1\in K$, let $\{x_n\}_{n\ge 1}$ be iteratively defined by $(1.4).$ Suppose that for any $x^*\in \displaystyle \bigcap\limits_{l=1}^{p}F(T_l)$, there exists a strictly increasing function $\Psi:[0,+\infty)\to [0,+\infty),\;\Psi(0)=0$ such that 
 $$(*)\;\;\;\;\Big<T_l^nx-x^*,j(x-x^*)\Big>\le k_n\|x_n-x^*\|^2-\Psi(\|x-x^*\|), \;{\rm for\; all\;}x\in K,\;l=1,2,...,p;\;p\ge 2.$$ Then $\{x_n\}_{n\ge 1}$ converges strongly to $x^*\in \displaystyle \bigcap\limits_{l=1}^{p}F(T_l).$
\end{theorem}

\begin{rmk}
There are a lot to say about this result but let us first and formost address the major issue arising from the proof of this theorem.\\

\noindent On page 37 of \cite{Rafiq}, immediately after inequality $(2.7)$, the author wrote: \\`` From the condition $(i)$ and $(2.7),$ we obtain $$\displaystyle \lim_{n\to\infty}\|y_n^1-x_{n+1}\|=0,$$ and {\sf the uniform continuity of $T_1$} leads to $$\displaystyle \lim_{n\to\infty}\|T_1^ny_n^1-T_1^nx_{n+1}\|=0."$$ This claim is, however, not true. To see this, we consider the following example:
\end{rmk}

\begin{example}
Let $\mathbb{R}$ denote the set of real numbers endowed with usual topology. Define $T:\mathbb{R}\to\mathbb{R}$ by $Tx=2x\;\forall\;x\in \mathbb{R},$ then $$|Tx-Ty|=2|x-y|\;\;\forall\;x,y\in \mathbb{R}.$$ This implies that $T$ is a Lipschitz mapping with Lipschitz constant $L=2.$ Thus, $T$ is uniformly continuous since every Lipschitz map is uniformly continuous. Now, suppose $y_n^1=1+\frac{1}{n}$ and $x_{n+1}=1-\frac{1}{n}$ for all $n\ge 1,$ then $$|y_n^1-x_{n+1}|=\Big|\Big(1+\frac{1}{n}\Big)-\Big(1-\frac{1}{n}\Big)\Big|=\frac{2}{n}\to 0\;{\rm as\;}n\to\infty.$$ We now show that 
$$\displaystyle \lim_{n\to\infty}|T^ny_n^1-T^nx_{n+1}|\ne 0.$$ Observe that 
\begin{eqnarray*}
Ty_n^1&=&2y_n^1=2\Big(1+\frac{1}{n}\Big)=2+\frac{2}{n}\\
 T^2y_n^1&=&T(Ty_n^1)=2\Big(2+\frac{2}{n}\Big) =2^2+\frac{2^2}{n}\\
T^3y_n^1&=&T(T^2y_n^1)=2\Big(2^2+\frac{2^2}{n}\Big)=2^3+\frac{2^3}{n}\\
&\vdots&\\
T^ny_n^1&=&2^n+\frac{2^n}{n} \;{\rm for\;all\;}n\ge 1.
\end{eqnarray*} Similar computation gives $$T^n x_{n+1}=2^n-\frac{2^n}{n}\;{\rm for\;all}\;n\ge 1.$$ Thus,
$$|T^ny_n^1-T^n x_{n+1}|=\Big|\Big(2^n+\frac{2^n}{n}\Big)-\Big(2^n-\frac{2^n}{n}\Big)\Big|=\frac{2^{n+1}}{n}\;\forall\;n\ge 1.$$
It is easy to see (using mathematical induction) that $2^{n+1}\ge n\;\forall\;n\ge 1$. So, $$|T^ny_n^1-T^n x_{n+1}|=\frac{2^{n+1}}{n}\ge 1\;\forall\;n\ge 1.$$ Hence, $$\displaystyle \lim_{n\to\infty}|T^ny_n^1-T^nx_{n+1}|\ne 0.$$ This contradicts the claim of A. Rafiq \cite{Rafiq}.\\
\end{example}

\noindent To correct the error in the result of A. Rafiq, we shall rather assume that $T_1$ is uniformly $L$-Lipschitzian so that $$d_n=M\|T_1^ny_n^1-T_1^n x_{n+1}\|\le ML\|y_n^1-x_{n+1}\|\to 0\;{\rm as\;}n\to\infty.$$ The rest of the result follows as in \cite{Rafiq}.

\begin{rmk}
In as much as the error in the proof of Theorem 5 of \cite{Rafiq} has been pointed out and corrected, it is not clear what the author really want to achieve by constructing such a complicated scheme given by $(1.4)$. If a clear study of the proof of \cite{Rafiq} is made, one will easily observe that the mappings $T_l, \;3\le l\le p$ played no role at all. This suggests that the scheme will only make sense if only two operators $T_1$ and $T_2$ are considered. Besides, it is not specified in Theorem 5 of \cite{Rafiq} which of the operators the sequence $\{k_n\}_{n\ge 1}$ is associated with. Meanwhile, condition $(*)$ gaurantees that the fixed point $x^*$ of these operators is unique. This thus reduces the entire problem to what has been studied in \cite{Chang} and \cite{Ofoedu}. We note that the result of Chidume and Chidume \cite{Chidume} and Ofoedu \cite{Ofoedu} remain correct if it were further assumed that the mapping $\phi:[0,+\infty)\to [0,+\infty)$ in thier results is onto.
\end{rmk}

\end{document}